\documentclass[11pt]{article}
  
  \usepackage[utf8]{inputenc}
  \usepackage[T1]{fontenc}
  \usepackage{geometry}
\usepackage{xcolor}
\usepackage[normalem]{ulem}
  \usepackage[frenchb,english]{babel}

\usepackage{chngcntr}
\counterwithin{subsection}{section}

  \usepackage{times}
  \title{Regards amis sur Claude Chevalley}
   \author{D. Couty}

\begin{document}
\selectlanguage{french}
  \maketitle

\stepcounter{section} 
 Nous nous proposons de suivre l'itinéraire  de Claude Chevalley lors des vingt dernières années de sa vie,  à travers les mots de Jacques Roubaud, Denis Guedj et Alexander\footnote{Selon Pierre Cartier, Grothendieck tenait beaucoup à cette orthographe \cite{Car00}.} Grothendieck,  nous replaçant dans  le contexte de leurs témoignages emplis d'amitié. Nous nous éloignons ainsi de l'image bien connue\footnote{Voir par exemple \cite{Dieu86} et \cite{Dieu87}.} de Claude Chevalley reposant sur ses rencontres essentielles, mathématiques et amicales\footnote{ Parlant d'amitié pour Claude Chevalley, il y eut bien sûr dès la rue d'Ulm ses rencontres avec Jacques Herbrand et Albert Lautman et nous ne pouvons oublier d'évoquer les  beaux textes, sobres et émouvants de l'ami japonais Shokichi Iyanaga \cite{Iya96}.}, au sein du groupe Bourbaki.

Nous terminerons ce parcours  en nous intéressant au devenir des nombreux écrits de Claude Chevalley que de grands mathématiciens ont  souhaité faire vivre après sa disparition.

\section*{ Le milieu des années soixante avec Jacques Roubaud}
\stepcounter{section}

 Jacques Roubaud  consacre quelques lignes intenses à Claude Chevalley dans \emph{`le grand incendie de londres', Branche 3  Mathématique :}  

 Au moment où il se place,  Chevalley est  professeur à la Faculté des Sciences de Paris.
\begin{quotation}\og L'année universitaire 64-65 le professeur Chevalley consacre son séminaire à la question de la `descente'. Il s'agit, dans son esprit, de donner un cadre catégorique épuré à la notion de `descente fidèlement plate' due à Grothendieck [...] J'ai admiré Chevalley énormément. Premièrement avant de le connaître, comme l'un des fondateurs de Bourbaki, particulièrement apprécié par moi parce que algébriste, et  parce que j'ai lu la rédaction, par mon ami Pierre Lusson, de son cours de 1958 à l'institut Henri-Poincaré, sur les formes quadratiques\footnote{Formes quadratiques sur un corps quelconque : 1er semestre 1955-1956. Cours Chevalley, Claude (1909-1984) / Faculté des sciences de Paris, centre de polycopie / 1955, site de l'idRef.}. Je l'ai apprécié aussi pour sa contribution à la fameuse théorie du corps de classe\footnote{La théorie du corps de classes est au centre de l'intérêt mathématique de Chevalley pendant la première partie de sa carrière.}[sic], roman aux deux héroïnes mathématiques mystérieuses, mesdemoiselles Adèle et Idèle. Je l'ai admiré deuxièmement pendant l'année  du séminaire sur la descente, où j'ai fait sa connaissance [...] Je l'ai admiré ensuite, troisièmement, comme homme remarquable. Je l'admire toujours, quatrièmement.\fg{} \cite{Rou09} \end{quotation}
Revenons sur les quatre points mis en évidence par Jacques Roubaud.

 \textbf{Premièrement avant de le connaître \dots{}}

 La découverte de Bourbaki par Jacques Roubaud  est marquée par l'empreinte de Pierre Lusson\footnote{D'après \cite{Rou09} pour l'ensemble des deux premiers  paragraphes.}. Pendant le début de l'année universitaire 1954 -1955, dans l'amphi Hermite de l'Institut Henri Poincaré, sans se connaître,  ils  assistent tous les deux au cours de Calcul différentiel et Intégral de \og M. G(ustave) Choquet, professeur\fg.  A l'automne 1954, Henri Cartan avait proposé Gustave Choquet pour  prendre  la suite de Georges Valiron\footnote{Georges Valiron (1884-1955) enseignait, d'après Jacques Roubaud, suivant le \emph{Cours d'Analyse mathématique} de Goursat dont le premier volume parut en 1902.} qui, malade, ne pouvait plus  assurer cet enseignement, 

 Choquet  en modifia résolument le contenu et l’orientation, introduisant la construction des nombres réels, les espaces topologiques, les espaces de Hilbert, au grand désarroi d’ailleurs des redoublants effarés par la métamorphose du sujet \footnote{ D'après la \emph{Notice nécrologique de Gustave Choquet lue par Michel Talagrand} lors de la séance publique de l'Académie des sciences du 2 octobre 2007.}.  Selon Jacques Roubaud, cet hiver là, l'auditoire était, sinon attentif, du moins particulièrement silencieux: \og ces silences avaient une densité et une tonalité particulières. Ils n'étaient indice ni d'émotion, ni d'enchantement, ni seulement de concentration appliquée. Ils marquaient avant tout la perplexité, ou même  la stupéfaction. Je partageais cette stupéfaction\fg{}. Pourtant, un jour, alors que nul jamais  n'interrompait  le \og  courant de la parole professorale magique\fg, une voix s'éleva faisant entendre un  \og mais\fg{}  début d'une intervention s'adressant au professeur Choquet. Choquet répondit et enchaîna, comme si tout cela était naturel.
Cependant, à l'issue du cours, cet épisode inhabituel suscita une discussion entre étudiants. A côté de celui qui avait  créé l’événement surgit un autre protagoniste qui semblait détenir une clé. Pour lui,  ce qui  sous-tendait le cours de Choquet,  ce qui au delà redonnait à \og \underline{La} Mathématique\fg{} son unité et son élan avait un nom : Bourbaki.  

C'est ainsi que  Bourbaki et Pierre Lusson firent  irruption dans la vie de Jacques Roubaud.  Dès les premiers mois de 1955, Jacques Roubaud passa  dans un \og grand calme studieux\fg{} ses soirées à la bibliothèque de la Sorbonne s'imprégnant du livre \emph{Topologie générale} de Bourbaki. Dans le même temps \og Trois de ces étudiants, et trois seulement, sont devenus alors et sont longtemps  restés\fg{} ses amis. Pierre Lusson est l'un deux. Amis et complices pendant de longues années, ils commencèrent presque ensemble leur carrière dans l'enseignement supérieur à  Rennes, Pierre Lusson en 1957 et Jacques Roubaud en 1958.\begin{quotation}\og [Jacques Roubaud,  est] 'entré dans la carrière'\dots{} universitaire à l'automne 1958, comme assistant délégué de mathématiques auprès de la faculté des sciences de l'université de Rennes.\fg {} \cite{Rou09}\end{quotation} En dehors des cours et des séances d'exercices consacrées aux étudiants,  Jacques Roubaud fit \og beaucoup d'algèbre\fg{} - c'est bien l'algèbre qui l'\og attirait vers la mathématique\fg. Dans le droit fil du bourbakisme, il s'approprie entre autres \emph{ Formes quadratiques sur un corps quelconque} - grâce  à \og{ un exemplaire  ronéoté et dédicacé par le rédacteur\fg{}  Pierre Lusson - premier cours de Chevalley à l'IHP  \og revenu de son exil aux USA\fg.

Quant aux héroïnes Adèle et Idèle, ainsi nommées en 1968 dans le faire part de décès de Nicolas Bourbaki\footnote{Qui fut \og distribué en 1968 à la sortie de quelques amphis par un mathématicien contestataire.\fg \cite{Chou95}}, la première des deux  notions  introduites est celle d'\emph{idèle}\footnote{La notion d'\emph{adèle} est due à André Weil \cite{Dieu99}.} - au masculin - seule due à Chevalley et dont le premier nom fut \emph{élément idéal}.   Chevalley définit l\emph{'élément idéal} en 1936 \cite{Che36},  le nom \emph{idèle} apparaissant officiellement en 1940 \cite{Che40}.   D'après Iyanaga \cite{Iya06},  la dénomination \emph{idèle}\footnote{ Peut-être les premières lettres de \emph{ideal Element}?} lui  aurait été suggérée par  Helmut Hasse : effectivement, après avoir présenté la notion d'\emph{élément idéal} dès 1935 dans une lettre adressée à Hasse, le nom \emph{idèle} apparaît   en 1937 dans un autre courrier de Chevalley à Hasse \cite{Roq13}.

\textbf{ Je l'ai admiré deuxièmement \dots{}}

La rencontre entre Jacques Roubaud et Claude Chevalley eut lieu  pendant l'hiver 1964-65 lors de son séminaire pour lequel : \begin{quotation}\og  Il sollicite la participation d'Adrien Douady, de Michèle Vergne\footnote{Elle a soutenu sa thèse \emph{Variétés des algèbres de Lie nilpotentes}  en juin 66, avec Chevalley pour directeur de thèse, avant de passer son doctorat d'état en mai 1971.  Elle a obtenu les prix Bordin et  Ampère de l'Académie des Sciences. Elle est directeur de recherche émérite au CNRS et membre de l'Académie des sciences.}, de Jean Bénabou\footnote{Spécialiste des catégories, il a préparé sa thèse  sous la direction d'Ehresmann. Grâce à l'impulsion de Chevalley, il finit par la passer en l'absence d'Ehresmann, parti pour une année sabbatique.}
 et de moi-même, accessoirement. Nous acceptons, Jean et moi, avec empressement.\fg{}  \cite{Rou09}\end{quotation} Jean Bénabou, nommé chargé d'enseignement au département de mathématiques de Rennes à la rentrée 63, est un passionné de mathématiques en quête à ce moment là d'une compréhension profonde de \og madame CATÉGORIE\fg{} \cite{Rou09}. Il entraîne Jacques Roubaud dans une \og exploration du monde catégorique\fg{} lui présentant, lors de fréquents échanges, le développement de ses idées mathématiques. Cela conduira  Roubaud  au séminaire Chevalley.
 \begin{quotation}\og J'étais là. J'étais pénétré de l'honneur qui m'était fait de participer à ce séminaire. Je n'étais rien, en particulier mathématiquement rien. Même pas normalien. Jean s'était porté garant de moi. De toute façon, le professeur Chevalley n'accordait qu'une attention très distraite aux hiérarchies.\fg{}  \cite{Rou09} \end{quotation}

Par ailleurs, Jean-Paul Benzécri\footnote{ Il soutint sa thèse en 1960  sur les variétés affines sous la direction d'Henri Cartan. Il est le fondateur de l'école française de l'analyse des données, développant des outils statistiques, notamment l'analyse factorielle des correspondances qui permet de traiter de grandes masses de données.},
nommé à la faculté de Rennes,  y fait au début des années 60 un cours de linguistique mathématique,  auquel assiste Jacques Roubaud. 

C'est à la rencontre de cette \og théorie mathématique de la syntaxe des langues naturelles\fg{} de Benzécri et des catégories de Bénabou que  Roubaud élabore tout au long de l'année 1965 un travail très personnel de recherche mathématique. 
 \begin{quotation} \og Mais qu'en faire? Encouragé par Jean, j'allais en tremblant présenter mes principaux objets et résultats à Chevalley, qui m'accueillit avec beaucoup de gentillesse, me conseilla de condenser le tout de façon à en faire deux notes aux  \uline{Comptes rendus de l'Académie des sciences de Paris} qu'il se chargea de `présenter' [...] \fg {} \cite{Rou09}\end{quotation}
Pendant  \og la rédaction de ce qui n'était encore qu'une rédaction sans finalité claire\fg, en novembre et décembre 1966, Roubaud expose ses travaux à la VIe section de l’École pratique des hautes études\footnote{Actuellement EHESS.}. C'est ensuite qu'il  demande à  Benzécri d'accepter son travail comme thèse, thèse soutenue le 17 février 1967.

Si pendant l'année universitaire 1964-65  lors de la préparation de ses exposés pour le séminaire, Roubaud se rend chez  Chevalley, c'est une autre raison qui va l'y conduire par la suite: \begin{quotation}\og  [...] j'eus de nombreuses occasions d'aller chez lui rue de Prony. Et je continuai à m'y rendre ensuite, jusqu'à la fin du séminaire et ensuite, parce qu'il me persuada de me mettre au jeu de go.\fg{}  \cite{Rou09}\end{quotation}

 \begin{quotation}\og \dots{} il se trouve qu'il avait appris à jouer au go au Japon et puis, à Paris, il ne trouvait pas de joueur [...] J'ai joué au go avec lui [...] et puis à un certain moment, on s'est dit, Pierre Lusson et moi-même, ça serait quand même bien de créer des circonstances telles que Chevalley puisse avoir des joueurs. Et donc, on a eu plein d'ambition, on s'est dit: "On va faire un traité de go, et à ce moment là plein de gens se mettront à jouer au go". \fg{}  \cite{Dug14} \end{quotation}

C'est ainsi que vint au monde le  \emph{Petit traité invitant à la découverte de l'art subtil du Go}  \cite{Lus69} qui fut écrit,  non sans humour, dans le  jardin fleuri du moulin d'Andé \footnote{Le moulin d'Andé, situé dans un boucle de la Seine est un très beau lieu, qui accueille depuis 1962 des artistes en résidence qu’il s’agisse d’écrivains, de cinéastes, comédiens ou musiciens. Il servit de décor au film de Truffaut \emph{Jules et Jim}. C'est là que résidait Georges Perec à l'époque de la rédaction du \emph{Petit Traité invitant \dots}}, par Pierre Lusson, Georges Perec et Jacques Roubaud.

\textbf{ Je l'ai admiré ensuite, troisièmement \dots{}}   Comment l'auteur du livre\emph{ `le grand incendie de londres'} a-t-il été touché par les qualités humaines de Claude Chevalley au point de nous dire qu'il est un homme remarquable ?  C'est en revenant à la jeunesse de ses parents que  nous connaîtrons la réponse de Jacques Roubaud.

 A la fin des années vingt, les parents de Jacques Roubaud,  Lucien Roubaud et Suzette Molino, lui, le philosophe, et  elle, l’angliciste, se rencontrent à l'ENS de la rue d'Ulm. Ils sont tous les deux de la promotion 1927 et  ont déjà passé auparavant une année de khâgne dans la même classe à Marseille, où: \og On ne se disait pas un mot, d'ailleurs: à ce moment là, les garçons appelaient les filles "mademoiselle". On ne se parlait pas. \fg{} \cite{Rou15} Ils étaient issus tous deux de milieux modestes. Pour l'un comme pour l'autre, l'entrée à l'ENS était un saut dans un autre milieu que le leur. Lucien Roubaud se souvient  qu'au moment où  certains de ses professeurs le poussaient vers les études, quelqu'un avait dit à son oncle, qui était son tuteur : \og La Rue d'Ulm?  Mais ce n'est pas pour des gens de votre origine \dots{}\fg{} \cite{Rou15}   Quant à  Suzette Molino, elle est, avec Clémence Ramnoux  et Simone Pétremont, l'une des trois premières jeunes filles reçues en lettres à l'ENS rue d'Ulm  \cite{Eft03}. Simone Weil y entrera l'année suivante.

Jacques Roubaud  dit de ses parents: 
 \begin{quotation}\og  [...] ils n’avaient que peu de rapports avec les élèves scientifiques qu’ils trouvaient généralement prétentieux et méprisants envers les littéraires, à l’exception, disait mon père, de Claude Chevalley\footnote{Extrait de mails échangés avec J.Roubaud, ainsi que pour les citations du \emph{quatrièmement}.}.\fg\end{quotation}

   La qualité du regard de Claude Chevalley évoquée au dessus, Jacques Roubaud  l'appelle sa \og{} grande modestie\fg. Elle se  trouvera confirmée pour lui, bien des  années plus tard, lors de ses rencontres avec le \og professeur Chevalley\fg. C'est ce qui sous-tend le `troisièmement'.

\textbf{Je l'admire toujours, quatrièmement.}

 A quelle source Jacques Roubaud  a-t-il puisé la force de cette affirmation, venant après les trois précédentes? \begin{quotation}\og Parce que jusqu’à sa mort il est resté fidèle à lui-même,  même quand les positions qui furent les siennes après 1968 ont été jugées défavorablement par le milieu mathématique.\fg\end{quotation}

 Cette rupture avec le milieu mathématique\footnote{Rappelons que  son  exclusion du prix  de l'Académie des sciences destiné au groupe 
 Bourbaki est encore très proche en 1968  - elle le touche directement de décembre 1966   à l'automne 1967 \cite{Cou19}.} dont parle Jacques Roubaud n'a pas été éludée par le milieu bourbakiste.  Elle est évoquée par Jean Dieudonné dans son article de 1986  \cite{Dieu99}:\begin{quotation}\og Jusqu'à la fin de sa vie, il n'a cessé de s'enflammer pour les victimes d'injustices [...] avec le plus parfait dédain des inimitiés qu'il aurait pu ainsi encourir. Il n'avait d'ailleurs que mépris pour les "honneurs", et refusait ceux auxquels il eût pu légitimement prétendre.\fg \end{quotation}
Ceci nous amène tout naturellement à la période qui commence en 1968  et aux  textes de Denis Guedj et d'Alexander Grothendieck.\\

\section*{ L'aventure de 1968 avec Denis Guedj}
  Dans les années 80,  Denis Guedj a réalisé des interviews de  Claude Chevalley dont, pour le moment, seuls quelques extraits sont accessibles et servent de trame à ce paragraphe\footnote{ Voir \cite{Gue85}, \cite{Gue04}, \cite{Gue04b}, \cite{Gue06}.}. 

\textbf{Printemps 1968}

En 1968, le travail de recherche que mène Denis Guedj   sous la  direction de Jean- Paul Benzécri concerne la théorie des grammaires formelles. Le Comité de grève s'installe alors dans le  batiment où il travaille \cite{Pe09}. C'est dans ce contexte qu'il  fait  connaissance avec Claude Chevalley :

\begin{quotation}\og Dans la faculté des sciences de Jussieu encore endormie\footnote{Depuis le milieu des années 50, les bâtiments de la Sorbonne  étant devenus trop petits, la faculté des sciences s'est implantée dans deux bâtiments  situés sur le site de l'ancienne Halle aux vins, quai Saint-Bernard et  rue Cuvier. Ce sont les premiers bâtiments du campus de Jussieu \cite{Del}.}, je passais devant le local du Comité de Grève. Un bruit. Dans la salle vide, avant les interminables réunions de la journée, un homme balayait consciencieusement le sol recouvert de mégots et de papiers: Claude Chevalley \dots{} Le Comité de grève avait pris place dans le laboratoire de linguistique mathématique un bâtiment préfabriqué, là où aujourd'hui s'élève l'Institut du Monde Arabe.\fg{} \cite{Gue04b}\end{quotation}

 Dans ce monde là, Denis Guedj et Claude Chevalley se sont longuement côtoyés:
\begin{quotation}\og Claude Chevalley a été l'un  des trois professeurs de la faculté des sciences à s'engager totalement dans l'aventure jusqu'à la fin, occupant les locaux avec les étudiants quai Saint-Bernard [...] et  y dormant fréquemment. C'est là que je l'ai rencontré.\fg{}\cite{Paj11} \end{quotation}

 \begin{quotation}\og  Nous avions pris possession de cet univers qui jusqu'alors n'avait été qu'un lieu d'études et de connaissances, et qui, dans la douceur de ce mois de mai, était devenu un lieu de vie, d'une vie merveilleusement grisante. La fac était à nous. La nuit, nous marchions dans les allées encore? longées de grands arbres, pénétrions dans les amphis vides, dormions à la belle étoile. Inutile de dire qu'à la rentrée, en automne 1968, il nous fut impossible de trouver notre place dans ces espaces déshabillés d'où la magie s'était retirée. \fg{}  \cite{Gue04b} \end{quotation}

 Denis Guedj est alors au CNRS et enseigne le traitement du signal en ${3}^{e}$ cycle. Il  raconte:

\begin{quotation}\og A la rentrée de septembre, impossible d'imaginer que je pourrais demeurer dans cette faculté redevenue "normale" alors que j'y avais vécu les moments parmi  les plus intenses de ma vie. Chevalley  était habité du même sentiment. Avec lui nous entreprîmes d'entretenir ailleurs cet esprit de liberté qui nous avait galvanisés.\fg{}  \cite{Paj11}\end{quotation}
 Claude Chevalley s'était proposé auprès des instances universitaires, c'est à dire auprès d'Edgar Faure pour partir dans les nouvelles créations\footnote{D'après Denis Guedj - De la Ronde au Zéro, Film couleur réalisé par Yolande Robveille et Patrice Besnard / 2010.}. Denis Guedj prend alors une décision :
\og J'ai démissionné du CNRS et j'ai été nommé maître assistant à l'université.\fg

 \textbf{L'aventure de Vincennes}

C'est à l'Université de Vincennes que les deux hommes passeront ensuite de longues années côte à côte\footnote{Ils partageront aussi l'expérience du mouvement \emph{Survivre} que nous présenterons au paragraphe suivant.} jusqu'au départ à la retraite de Claude Chevalley en 1978.

 Un projet d'université nouvelle s'était  concrétisé  sous la forme du Centre expérimental de Vincennes. Il était  le fruit du bouillonnement intellectuel du printemps 1968.

\og L'université de Vincennes n'est pas sortie, comme par miracle, du néant\fg{}  écrit  Raymond Las Vergnas, professeur de langue et littérature anglaises\footnote{ Pendant le Front populaire, Las Vergnas fut  chef de cabinet adjoint de Jean Zay,  ministre de l’Éducation nationale et des Beaux-arts - voir le fonds Jean Zay sur le site des Archives nationales.},  doyen de  la faculté des lettres de Paris en  1968. 
Le 5 août 1968, il présenta au nouveau ministre de l'éducation nationale, Edgar Faure,  son projet d'université nouvelle:
 \begin{quotation}\og L'université nouvelle, ai-je dit au ministre, serait conçue selon un principe fondamental de participation que l'on retrouverait à tous les niveaux, tant de la gestion que de la pédagogie [...] Quelques jours plus tard, M. Edgar Faure m'annonçait qu'il me donnait carte blanche pour créer le nouveau centre expérimental [...] Il fallut attendre le 10 décembre pour que le \emph{Journal officiel} se décide à rendre publique sa création alors que les bâtiments, construits, eux, en un temps record, étaient sortis de terre depuis des semaines.\fg{}  \cite{Co79}\end{quotation}

Quelques idées président, selon lui, à la création de cette université: adaptation de l'enseignement à l'évolution du monde par une augmentation de la part des exercices pratiques, contrôle continu des connaissances, personnalisation des études, remplacement des certificats de licence par des unités de valeur,  assouplissement des conditions d'entrée pour les non-bacheliers,  recrutement des enseignants selon des critères échappant \og aux filières académiques rituelles\footnote{Il fut en effet décidé de créer un \og noyau cooptant\fg, composé d'une trentaine d'enseignants de toutes disciplines, reconnus pour leurs compétences et animés d'une vocation novatrice. Ces enseignants eurent la charge de recruter sur dossiers concurrentiels environ 200 autres enseignants \cite{Co79}.}\fg.
\begin{quotation}<<Une nouvelle université venait de s'ouvrir dans le bois de Vincennes. Nous avons joyeusement émigré. Claude y est resté jusqu'à sa retraite [...] Nous avons créé le département de mathématiques [...] Pour Claude, ce furent "des années de bonheur, les années les plus heureuses de ma vie".\fg{} \cite{Gue04b} \end{quotation}

Dans cette université, \og expérimentation très libre et très ouverte\fg{} \cite{Co79},  prévue pour s'ouvrir aux salariés, pour faire vivre la  pluridisciplinarité, des cours d'alphabétisation avaient été proposés aux travailleurs de la faculté.
 \begin{quotation}\og Durant plusieurs semaines, avec le sérieux qui le caractérisait, il [Chevalley] a donné au seul travailleur inscrit, un jeune malien employé au service de nettoyage, un cours particulier sur l'addition, la multiplication. Voilà comment, aux dires de nombre de ses collègues, Chevalley perdait son temps\footnote{ Il est vrai que, par exemple, au même moment, Dieudonné  est  à l'Académie des sciences où il a été élu  le 24 juin 1968.} [...] \fg{} \cite{Gue04b} \end{quotation}

 Denis Guedj et Chevalley  quant à eux, se voient très souvent:\begin{quotation}\og Nous avons beaucoup réfléchi et appris ensemble; il est celui qui m'a le plus profondément marqué et pour lequel j'ai une affection profonde. Un père-frère, un complice.\fg{} \cite{Paj11}\end{quotation}

\begin{quotation}\og Depuis cette époque nous sommes restés très proches. J'ai pour lui une infinie tendresse et un profond respect.\fg{} \cite{Gue04b}\end{quotation}

 Denis Guedj, questionné par le rapport entre sciences et société,  se passionne pour la philosophie, l'histoire, l'histoire des sciences. 
\begin{quotation}\og Une grande partie de ce trajet, je l'ai accomplie aux côtés de Claude Chevalley [...] ce fut, et cela reste, la rencontre capitale. Un cadeau qu'il me faisait, un cadeau que je me faisais. C'est avec lui que j'ai éprouvé ce que le contenu même du savoir mathématique pouvait comporter d'émotion. Étonnant quand on a à l'esprit l'image de Chevalley cofondateur du groupe Bourbaki.\fg{} \cite{Paj11}\end{quotation}

Il nous  livre ce  portrait:

\begin{quotation}\og Chevalley, c'était une apparence fragile, une pensée ferme, ouverte au questionnement [...] Peu disposé à accepter les injustices et la violence faite aux faibles, il était tout sauf modéré. Il se situait à l'intersection de quatre directions. Les mathématiques, la philosophie, l'engagement politique, la foi\footnote{Même si  Chevalley, chrétien engagé, s'est éloigné du protestantisme dans les années 70 \cite{Gue04b}.} l'ont mobilisé et  l'ont constitué [...] Liberté était son maître mot  [...] \fg{}
\cite{Gue06}\end{quotation}

\section*{Le tournant des année soixante-dix par Alexander Grothendieck}

\textbf{Les premiers moments au sein de Bourbaki}

 Pour Alexander Grothendieck\footnote{ Pour approcher  la personnalité et l’œuvre  de Grothendieck,  on peut profiter du regard de Pierre Cartier (\cite{Car00}, \cite{Car09}).}, le nom de Chevalley est d'abord associé à celui des autres membres de Bourbaki.\

Alexander Grothendieck  arrive à vingt ans, en 1948,  à Paris, avec dans sa \og{maigre valise une Licence es Sciences de l'Université de Montpellier}\fg. Sur la recommandation d'un de ses professeurs de Montpellier, il rencontre Henri Cartan, dont il va suivre le séminaire.\

 Décrivant \og l'étranger bienvenu\fg{}
 qu'il fut lui-même au sein de Bourbaki , il dit:\ 
\begin{quotation}\og Au Séminaire Cartan il y avait aussi des apparitions périodiques de Chevalley, de Weil, et les jours des Séminaires Bourbaki (réunissant  une 
petite vingtaine ou trentaine à tout casser, de participants et auditeurs), on y voyait débarquer, tel un groupe de copains un peu bruyants, les autres membres de ce fameux gang Bourbaki : Dieudonné, Schwartz, Godement, Delsarte.\fg{} \cite{Gro86} \end{quotation}

 Claude Chevalley, alors aux Etats-Unis, est rentré à Paris pour l'année universitaire 1948-49 grâce à une bourse Guggenheim \cite{Dieu99}.

 Il expose par exemple en décembre 1948 au séminaire Bourbaki sur \emph{L'hypothèse de Riemann pour les corps de fonctions algébriques de caractéristique p} \footnote{Chevalley codirigera le séminaire Cartan à son retour en France, pendant l'année universitaire 1955-56 et fait plusieurs  exposés sur les schémas - d'après \emph{Tome 8,  séminaire Henri Cartan}.}.

 A ce moment-là, pour Grothendieck,  Chevalley est un  Bourbaki parmi les autres.

 \textbf{ Survivre \dots{} et Vivre d'après Alexander Grothendieck}  

 L'essentiel de la rencontre  entre Grothendieck et  Chevalley se passe en dehors de Bourbaki et laissera une trace profonde chez Grothendieck.
Dans  \emph{ Récoltes et Semailles}, le nom de Chevalley revient à de multiples reprises sous sa plume.
Dès la page 29, il évoque:
\begin{quotation}\og  [...] Claude Chevalley, le collègue et ami à qui est dédiée la partie centrale de Récoltes et Semailles [...] En plusieurs endroits de la réflexion, je parle  de lui, et du rôle qui fût le sien dans mon itinéraire.\fg
\end{quotation}

C'est au début des années 70 que  Grothendieck et Chevalley  vont faire plus ample connaissance  au sein du mouvement \emph{Survivre \dots{} et Vivre} \cite{Gro86}. 

\emph {Survivre \dots{} et Vivre} est le nom d’une revue\footnote{Dont le nom de départ  \emph{Survivre}  est assorti dès le \no2 d'un sous-titre quelque peu  emphatique: \emph{Mouvement international et inter professionnel pour notre survie, fondé le 20.7.1970 à Montréal}. L'histoire de ce mouvement éphémère a été étudiée par Céline Pessis  \cite{Pe09}.} et d’un mouvement pacifiste et écologiste  créé au début des années 70 par quelques scientifiques. La revue est en fait un petit journal\footnote{1300 exemplaires pour le \no 6, 12 500 pour le \no12 de juin 1972 d'après \cite{Pe09}. Un certain nombre de ces journaux est consultable en ligne: http://www.grothendieckcircle.org/} d'une quarantaine de pages, dont seulement 19 numéros paraîtront, le mouvement se prolongeant après la disparition de la revue. Le ton,  la forme du journal, la facture des dessins\footnote{On peut y voir  par exemple  les premiers dessins de Didier Savard - \nos 8 et 12 -  qui  publie un peu plus tard ses dessins de presse pour \emph{Libération}, \emph{La gueule ouverte} \dots{} Il  devint par la suite un auteur reconnu d'albums de BD.} sont en phase avec  la mouvance anti-conformiste de l'après 68.

 Il est permis de penser que Claude Chevalley  a retrouvé dans ce mouvement certains des élans qui l'avaient porté vers le mouvement \emph{ l'Ordre nouveau} de sa jeunesse. En effet, dans le \no2/3 daté de septembre/octobre 1970, à la rubrique  \emph{Des adhérents se présentent}, il choisit de dire:\begin{quotation}\og J'ai participé durant la même période\footnote{Il s'agit de la période 1931-1937.} à un mouvement politique appelé "l'Ordre nouveau" (qui n'eut rien de commun avec les mouvements qui reprirent ce nom par la suite) dont la tendance dominante était le personnalisme teinté de certaines influences anarchisantes.\fg{} \cite{Che70}
\end{quotation}
 Bien que Grothendieck dise  de Chevalley qu'il \og{s'était joint  au groupe avec une conviction mitigée}\fg, il semble s'y être assez fortement impliqué car on peut lire, en première page, à partir du numéro de l'automne 1970 - et jusqu'à décembre 71 -  \og{Directeur de  publication (édition française):  Chevalley}\fg. Il fait aussi partie, avec  Grothendieck, du comité de rédaction  et tient des permanences à son domicile les lundis après-midi. Guedj  suivra Chevalley à \emph{Survivre} par amitié \cite{Pe09}.

 Après  \emph{l'Ordre nouveau}, après Bourbaki,   après  \emph{Vincennes}, le voila à nouveau engagé dans un mouvement collectif. Encore une fois, il fait partie d'un petit groupe qui a pour projet de faire émerger des idées nouvelles. 

 Au fil des numéros, Claude Chevalley signe plusieurs articles, dont des analyses de livres\footnote{Par exemple,  \emph{Une société sans école}, d'Ivan Illich à la rubrique \emph{Le livre du mois} du \no 12 de \emph{Survivre \dots{} et vivre}.}. En janvier 1971\footnote{La date indiquée dans ce numéro  est janvier 70, mais le numéro précédent étant daté de décembre 1970, il  s'agit visiblement d'une erreur de frappe.}, dans le \no 6 de \emph{Survivre},  Claude Chevalley  choisit de présenter \emph{The New Brahmins, Scientific Life in America} de  Spencer Klaw \cite{Kl68}. Le commentaire de Claude Chevalley sur cette étude de la situation des chercheurs scientifiques aux USA mérite d'être lue dans son intégralité, même si seulement quelques extraits sont cités ici:
\begin{quotation}\og Avec la liberté de décider lui-même de l'orientation de son travail, le chercheur perd naturellement tout intérêt pour ce qu'il fait. Faute d'avoir dès sa jeunesse réfléchi aux répercutions sociales possibles de son activité de chercheur, il est tout prêt à accepter que son travail n'ait d'autre lien avec la réalité que le profit que l'entreprise peut en tirer [...]

Le chercheur qui en est à ce point, quand il ne se désintéresse pas totalement de son activité professionnelle, s'oriente tout naturellement vers les activités administratives seules propres à lui conférer l'estime et la considération du milieu où il vit [...] il convient cependant de noter qu'il est souvent conscient d'être passé de l'autre côté de la barrière et d'avoir renoncé au sens qu'il avait donné à sa vie  en choisissant le métier de chercheur: conscience qui contribue encore à la rapidité de son évolution en lui faisant éviter tout rapport avec ceux de ses ex-collègues qui sont encore chercheurs.\fg\end{quotation}
  Claude Chevalley pense-t-il à ce moment-là aux \og{pontifes}\fg{} de sa jeunesse ou aux 
\og {mandarins}\fg{} des années 70 ? 

\textbf{L'éloignement de Bourbaki}

En 1986, quand Grothendieck s'interroge,  dans un  paragraphe   intitulé  \emph{Le mérite et le mépris}  sur le \og règne du mépris\fg{} qui a gagné, selon lui, le milieu mathématique, il nous fait partager  ce qui distingue sur ce point Chevalley du reste de la communauté mathématique.
\begin{quotation}\og C'est d'ailleurs Chevalley qui a été un des premiers, avec Denis Guedj que j'ai aussi connu par Survivre, à attirer mon attention sur cette idéologie-là (ils l'appelaient la "méritocratie" ou un nom comme ça), et ce qu'il y avait en elle de violence, de mépris. C'est à cause de ça, m'a dit Chevalley [...] qu'il ne supportait plus l'ambiance dans Bourbaki  et avait cessé d'y mettre les pieds.\fg{} \cite{Gro86}\end{quotation}
Il tient à réserver une place à part à  Chevalley, ajoutant bien plus loin au paragraphe \emph{Trois jalons - ou l’innocence}:

\begin{quotation}\og Il [un esprit de suffisance] a dû venir à pas de loup, au cours des ans, s’installant à demeure en les uns et en les autres, peu à peu, sans que personne
parmi nous (mis à part Chevalley seulement \dots{} ) ne s’en aperçoive.\fg   \end{quotation} 

 On peut rapprocher ces quelques lignes des propos de Chevalley décrivant le changement d'état d'esprit de Bourbaki, qu'il n'apprécie plus à partir d'un certain moment. \begin{quotation}\og L'esprit de canular, l'esprit de Jarry, a été rattrapé par une certitude  absolue de supériorité sur tous les autres mathématiciens.\fg{} \cite{Chou95} \end{quotation}
\begin{quotation}\og Ce qui m'a détaché principalement de Bourbaki, c'est de "le"  voir prendre des positions réactionnaires à l'Académie et dans les universités.\fg{} \cite{Gue06}\end{quotation}
Puis dans un article de 1982 qu'il cosigne avec M. Kasner, on peut lire:
\begin{quotation}\og La vie mathématique semble de plus en plus dominée par des clans qui s'y conduisent comme en pays conquis, s'emparant des positions-clés dans les institutions et périodiques mathématiques, appliquant en tous lieux le projet d'Armande: "Nul n'aura de l'esprit hors nous et  nos amis".\fg{} \end{quotation}
Et y incluant sévèrement Bourbaki :

\begin{quotation}\og [...]  histoire classique des jeunes révolutionnaires  à qui l'âge a permis de prendre le pouvoir.\fg{}  \cite{Che82} \end{quotation}

\textbf{Regards et derniers mots sur Claude Chevalley}

Au delà de l'histoire du groupe \emph{Survivre \dots{} et Vivre}, au delà de l'amertume  de Chevalley sur le milieu mathématique,  Grothendieck  nous livre une image de  Chevalley  \og au delà des façades de rigueur\fg{}  :

\begin{quotation}\og Il lui arrivait parfois  de parler de lui-même, juste quelques mots à l'occasion de ceci ou  cela, avec une simplicité déconcertante [...] Il parlait peu, et ce qu’il disait exprimait, non des idées qu’il aurait adoptées et faites siennes, mais une perception et une compréhension personnelle des choses [...]  Ce qu'il disait bousculait souvent des façons de voir qui m'étaient chères, et que pour cette raison je considérais comme "vraies" [...]

Je me rendais compte obscurément qu'il avait quelque chose à m'apprendre sur la liberté - sur la liberté intérieure.\fg{} \cite{Gro86}\end{quotation}

En juin 1984, Grothendieck pense être sur le point de terminer  \emph{Récoltes et Semailles}, prévoyant seulement la rédaction de deux ultimes notes.

\begin{quotation}\og Ça faisait des mois que je me voyais sur le point d’en terminer avec Récoltes et Semailles, frappé tiré broché et tout - et de monter à Paris dare dare pour lui apporter un exemplaire encore tout chaud ! S’il y avait une personne au monde dont j’étais sûr qu’elle lirait mon pavé avec un vrai intérêt, et avec plaisir souvent, c’était lui - et je n’étais pas sûr du tout s’il y en aurait un autre que lui ! 

 Dès les débuts de ma réflexion, je m'étais rendu compte que Chevalley m'avait apporté quelque chose, à un
moment crucial de mon itinéraire, quelque chose semé dans une effervescence, et qui avait germé en silence [...]

 De le rencontrer et de parler avec lui tant soit peu m'aurait permis sûrement de mieux appréhender cet ami que par le passé, et de mieux situer et cette parenté essentielle, et nos différences. S'il y avait, à part Pierre Deligne\footnote{ Pierre Deligne est un mathématicien belge, né en 1944. Élève de Grothendieck, il a soutenu sa thèse d’État en 1972. Il a reçu de nombreux prix mathématiques, en particulier la médaille Fields en 1978 pour sa résolution des conjectures de Weil et  le prix Abel en 2013  pour l'ensemble de ses travaux.}, une personne pour laquelle je ressentais une hâte de pouvoir lui remettre en mains propres le texte de Récoltes et Semailles, c'était bien Claude Chevalley. S'il y avait une personne dont le commentaire, espiègle ou sarcastique, aurait pour moi un poids particulier, c'était lui encore. \fg  \end{quotation}

Mais rien ne se passa comme prévu.

\begin{quotation}\og En ce jour-là de la première semaine de juillet, j'ai su que je n'aurais pas ce plaisir de lui apporter ce que j'avais de meilleur à offrir, ni celui d'entendre encore le son de sa voix.\fg{} \cite{Gro86}\end{quotation}
Claude Chevalley est mort le  28 juin 1984.

\section*{Après 1984, regards mathématiques sur son œuvre} 

 Le monde mathématique est présent autour de  Claude Chevalley disparu. Deux des fondateurs de Bourbaki tout d'abord :  Dieudonné envoie à Weil son  projet de note biographique\footnote{Voir le \emph{Dossier biographique} de Weil à l'Académie des Sciences.}  sur Chevalley. Weil répond par une relecture détaillée \cite{Wei85}, lui conseillant en particulier d'aborder les aspects familiaux de la vie de Claude Chevalley, son mariage avec Jacqueline,  son remariage avec  Sylvie, sans oublier bien sûr leur fille Catherine. Ainsi voit le jour en 1986  la notice biographique de Dieudonné  pour  l'Annuaire des Anciens  élèves de l’École Normale Supérieure \cite{Dieu99}. Puis vinrent  les articles de Dieudonné dans la\emph{ Vie des sciences} \cite{Dieu86},  Dieudonné et Tits dans le bulletin de l'AMS \cite{Dieu87}.

De nombreux mathématiciens vont ensuite  se  mobiliser pour publier  l’œuvre mathématique de Claude Chevalley. C'est une mission difficile.

 Nous savons  grâce à  Pierre Cartier et à Catherine Chevalley \cite{Che97} que Chevalley avait, dès 1982,  exprimé  trois souhaits  concernant  la publication de ses travaux.
 Il aurait aimé qu'elle puisse contenir ses écrits \emph{non- techniques}, épistémologiques ou politiques. Il aurait également souhaité que ses écrits mathématiques, dont certains  passages lui paraissaient insatisfaisants,  puissent être révisés par des notes critiques. Son troisième souhait était d'y inclure certains textes jamais publiés, mathématiques et non-mathématiques.  

Chevalley s'était ouvert à eux sur les raisons profondes de ces choix et il avait pleinement conscience que cela rendait ardue la publication de ses écrits, car il n'imaginait pas se  replonger dans ses anciens travaux, pas plus qu'imposer un tel pensum à  autrui.

  Cela rejoint les mots de l'ami japonais Iyanaga, auquel,  un an avant sa mort, Claude Chevalley s'était confié: \begin{quotation}\og Claude m'a dit que la Maison Birkhäuser lui avait proposé de publier ses œuvres, mais qu'il ne voulait pas laisser paraître tous ces papiers en photocopie parce qu'il n'en aimait pas certaines et que, d'un autre côté, il se trouvait des inédits qu'il aurait bien voulu faire publier \dots{}\fg{} \cite{Iya96}\end{quotation} 

C'est dans le même esprit  que Cartier rappelle  en 2005 : \begin{quotation}<< Chevalley lui-même avait insisté sur le fait qu'une publication de ses Œuvres supposerait une révision très soignée.>> \emph{Avertissement au lecteur} \cite{Che05} \end{quotation}

\stepcounter{section}
 \textbf{ Comité éditorial des \emph{Œuvres} de Chevalley}

Et pourtant,  après la mort  de Chevalley, on trouve  Cartier à la tête du comité éditorial  des \emph{Œuvres} de Chevalley. 

En 1997, le premier tome de cet \og ambitieux projet\fg {} soutenu par le CNRS  voit le jour. Il s'agit du \emph{Volume II}, intitulé  \emph{The Algebraic Theory of Spinors and Clifford Algebras}  \cite{Che97}. 

Les textes ont été revu avec le plus grand soin, conformément aux vœux de Chevalley. Pour ce volume,  Cartier 
 a pu compter  sur la collaboration de nombreux mathématiciens comme Michel Broué, Michel Enguehard et Jacques Tits, membres du \emph{Séminaire Chevalley}.  Jean-Pierre Serre, Armand Borel, Shokichi Iyanaga, Henri Cartan ont aussi participé à l'aventure. L'avant-propos est signé de Catherine Chevalley et de Pierre Cartier. L'ouvrage contient   également la critique par Dieudonné du livre  \emph{The Algebraic Theory of Spinors} de Chevalley de 1954 et une postface due à Jean-Pierre Bourguignon\footnote{Jean-Pierre Bourguignon y actualise  cet ouvrage  de Chevalley des années cinquante en le reliant à ses  applications en physique théorique. Cela n'est pas sans nous rappeler l'\emph{Appendice - Les progrès récents de la théorie des nombres} que Chevalley écrivit en 1936  lors de la publication posthume d'un mémoire d'Herbrand \cite{Her36}.}.

Dès cette période, le  découpage en six volumes des \emph{Œuvres} de Chevalley est établi:\\

I.  \emph{Class Field Theory}\hspace{0,5cm}     II.  \emph{Spinors} \hspace{0,5cm}  III.  \emph{Commutative Algebra and Algebraic Geometry}\hspace{0,5cm} IV. \emph{Algebraic Groups}\hspace{0,5cm} V  \emph{Epistemology and Politics}\hspace{0,5cm} VI.  \emph{Unpublished Material and Varia}\\

Nous avons quelques pistes sur le contenu des autres volumes:

Il est annoncé dans l'avant propos \cite{Che97}  que le \emph{Volume I}  contiendra les deux nécrologies de Tits et  Dieudonné, des lettres de Chevalley, la liste  même incomplète\footnote{Il y travaillait en 1984} de ses écrits non publiés et qu'il souhaitait faire paraître. Iyanaga serait chargé de l'introduction de ce \emph{Volume I}.

 Parmi les écrits du \emph{Volume V} prendraient place  des lettres de et à Jacques Herbrand et Emmy Noether.  
En 1998,  Cartier confirme \cite{Sen98}, qu'en tant qu'éditeur des \emph{Œuvres} de  Chevalley, il a prévu,  à la demande de Catherine Chevalley,  un volume rassemblant  les écrits non-mathématiques de  Chevalley, celle-ci  se chargeant de les réunir.

Une bibliographie complète des écrits de Chevalley  se trouverait dans le \emph{Volume VI} \cite{Che97}.

 La publication des \emph{Œuvres} paraît bien  lancée.

En 2005, Pierre Cartier, qui réussit le long et  délicat travail de réédition  du Séminaire  dirigé par Claude Chevalley  à l’École Normale Supérieure pendant les années universitaires 1956/57 et 1957/58,  s'interroge sur la place à donner à cet ouvrage:  \begin{quotation}\og Une fois prise la décision de publier le tout, nous hésitâmes un moment à en faire un des volumes des "Œuvres complètes" de Chevalley, vu la part importante des autres rédacteurs.\fg\cite{Che05}\end{quotation}
 Le titre fut finalement  \emph{Claude Chevalley} - \emph{Classification des Groupes Algébriques semi-simples - Collected Works}, Vol.~3.
\\

 \textbf{ Quelque chose d'inachevé}\\

Mais des blocages  stoppèrent le projet de publication  des  autres volumes des \emph{Œuvres}.

En 1999, parait une bibliographie très détaillée réalisée avec l'apport de documents transmis par Catherine Chevalley \cite{Dieu99}.  Cette liste\footnote{Elle ne recense pas tous les écrits de Chevalley. On pourrait en effet y ajouter les interviews réalisées par Denis Guedj au début des années 80,  certains textes du fonds Bourbaki - jamais publiés et auxquels tenait particulièrement  Chevalley - comme \emph{Introduction à la théorie des ensembles} et \emph{Géométrie élémentaire}  \cite{Che97} , certains textes de la revue \emph{Survivre}.}  réunit les références de nombreux textes, sous quatre  rubriques \emph{Mathematics},  \emph{Philosophy}, \emph{Publicity writing}, \emph{ Unpublished}: à côté des écrits mathématiques  s'y trouvent des textes philosophiques,  parfois co-écrits avec Arnaud Dandieu ou Alexandre Marc,  s'étalant de 1932 à 1962, les articles parus dans la revue \emph{Ordre Nouveau} entre 1933 et 1938,  puis d'autres écrits  plus tardifs, publiés dans diverses revues \emph{Réforme}, \emph{Dédales} \dots{} ainsi que des textes jamais publiés. Mais cet article paraît sans lien avec la publication des \emph{Œuvres}.

Le \emph{Volume I} ne parut pas:   Iyanaga apprit l'annulation du projet en 2005. C'est alors qu'il  décida de publier sous une autre forme le travail qu'il avait rédigé comme \emph{Introduction} à ce volume.  Ce  texte  \cite{Iya06}, publié en 2006, fut proposé au \emph{Japanese Journal of Mathematics} alors que Shokichi Iyanaga  était hospitalisé quelques mois avant son décès: ultime hommage de Shokichi  à Claude\footnote{Toshiyuki Kobayashi se chargea de cette démarche à la demande de Iyanaga.}.

Mais, malgré la volonté et les travail de nombreux mathématiciens, le projet de publication des \emph{Œuvres complètes} de Chevalley est à l'arrêt depuis 2005.
Le \emph{Volume II} et ce \emph{Volume III}  sont  à ce jour les seuls volumes  publiés : deux volumes sur les six prévus.  \\ 

 Cartier confie à ce propos  \og C'est une tragédie!\footnote{Lors d'un entretien privé en novembre 2016 à Pau.}\fg\\

 Par ailleurs, le devenir des entretiens réalisés par Guedj et  réunis dans \emph{Conversation avec Claude}\footnote{Que Guedj appelle aussi  \emph{Parler avec Claude Chevalley}.} reste incertain. En 1995, on peut lire :
 \begin{quotation}\og L'exploitation de ceux-ci, conservés par Catherine Chevalley, physicienne et historienne des sciences, devrait mieux faire comprendre l'évolution philosophique de l'un des grands mathématiciens de l'équipe [Bourbaki].\fg{} \cite{Chou95}\end{quotation}
 Mais en 2004, même si Guedj espère que \og \emph{Parler avec Claude} sera publié. \fg{} \cite{Gue04},  il dévoile les difficultés liées à la publication de ces interviews, ce qu'il confirme en 2008 :
 
\begin{quotation}\og Durant plus de dix années, deux longs après-midi chaque semaine, nous [Chevalley et lui] avons parlé ensemble. Parlé et pas écrit. Remords de ma part? J'ai dans les années 80 fait quelques interviews rassemblées sous le titre \emph{Conversation avec Claude Chevalley}. Mais ce recueil est inédit. Inédit parce que depuis sa mort, ses proches en ont empêché la diffusion. Peut-être faudra-t-il un jour passer outre cette interdiction.\fg{} \cite{Paj11} \end{quotation}

 En 2004,  seuls quelques extraits de ces entretiens  trouvèrent place  dans la revue \emph{Tangente}.\\

Quelque chose reste donc inachevé \dots{} \\

\section*{Et pour conclure}
Bouclant la boucle, revenons en 1964 et intéressons nous au regard que Chevalley pose sur Emil  Artin,  le maître de ses jeunes années:  selon Chevalley, on devrait définir  \og [...]  le fait d'être algébriste plus comme reflétant un certain tempérament intellectuel  que comme indiquant un sujet d'études privilégié.\fg{} Parlant d'Artin, il précise : \og [...] ce n'est pas seulement dans son activité mathématique que cette manière d'être qui était la sienne se manifestait : affirmation du primat de l'intellect par rapport à la passion, du conscient par rapport à l'inconscient,  de l'enquête méthodique par rapport aux éclairs de l'intuition, voilà une série de traits souvent associés et qui rapprocheraient par exemple Artin de Valéry ou de Mallarmé. Ce sont aussi sans doute ces  traits  et cette manière d'exister qui devaient rendre insupportable à Artin la tyrannie nazie, à laquelle il échappa dès que la chose lui fut matériellement possible.\fg{} \cite{Che64}

 Dans son regard sur Artin, Chevalley  affirme l'unité de pensée à laquelle il tient et qu'il a souhaité pour lui. La vision de l'entité constituée par sa vie et son œuvre,  qui se doit de rassembler mathématique, philosophie, épistémologie et politique, telle que l'a recueillie le comité éditorial  de ses \emph{Œuvres} reste fidèle à ses idéaux de jeunesse.\\

\end{document}